\input amstex
\documentstyle{amsppt}
\nologo
\advance\voffset  by -1.0cm
\NoBlackBoxes
\magnification=\magstep1
\hsize=17truecm
\vsize=23.2truecm
\voffset=0.5truecm
\input epsf
\document

\def\bR{{\Bbb R}}

\def \al {\alpha}

\def \part {\partial}

\leftline{\eightit ``I think you might do something better with the time'', Alice said, }

\leftline{\eightit  ``than wasting it in asking riddles that have no answer ...''}

\leftline{\eightit \quad - from L.~Carroll 'Alice's Adventures in Wonderland'}

\vskip 1.0cm

\topmatter
\title
A few riddles behind  Rolle's theorem
\endtitle

\author B.~Shapiro$^{\dag}$, M.~Shapiro$^\ddag$
\endauthor
\affil $^\dag$
      Department of Mathematics, University of Stockholm, S-10691, Sweden,
{\tt shapiro\@math.su.se}\\
$^\ddag$ Department of Mathematics, Michigan State University, East
Lansing, MI 48824-1027, USA, {\tt mshapiro\@math.msu.edu}\\
\endaffil
\leftheadtext {Shapiro \& Shapiro }
\rightheadtext {Riddles of Rolle's theorem}


\endtopmatter

\heading  Getting started \endheading

 First year undergraduates usually learn about the classical Rolle's
theorem saying that between two consecutive
zeros of a smooth function $f$ one can always find at least one zero
of its derivative $f^\prime$. In this paper we study a 
generalization of Rolle's theorem dealing with the zeros of higher
derivatives for a class of smooth functions which we call {\it $n$-nice}. 
For such functions we obtain some mysteriously looking additional 
inequalities governing the mutual arrangement of the zeros of different
derivatives.  The considered topic  is so classical that 
we find  it  impossible to be sure that our results are new. 
However, we were unable to trace anything similar in the literature 
    although such facts might have been known to at least 
  Chevalier Augustin Cauchy if not to  Sir Isaac Newton himself.

Consider a smooth function $f$ that on a certain interval has  $n$ distinct
real zeros  which we denote by
$x_{1}^{(0)} < x_{2}^{(0)} <\ldots <x_{n}^{(0)}$ . Then, by the
classical Rolle's theorem, $f'$ has at
least $(n-1)$ zeros, $f''$ has at least $(n-2)$ zeros, ... , $f^{(n-1)}$
has at least
one zero on $(x_{1}^{(0)},x_{n}^{(0)})$.

We call a smooth function $f$ having $n$ simple real zeros on some
interval an {\it $n$-nice function}  if for all $i=0,\ldots,n$ the
$i$-th derivative  $f^{(i)}$ has on the same interval exactly $(n-i)$
zeros denoted by $x^{(i)}_1<x_2^{(i)}<...<x_{n-i}^{(i)}.$ Note
that we require, in particular, that $f^{(n)}$ is nonvanishing!
Observe also that if $f$ is $n$-nice on $I$ then  for all $i<n$ 
its derivative $f^{(i)}$ is $(n-i)$-nice on the same interval. As a natural example of an $n$-nice
function on the whole $\bR$ one can take any polynomial of  degree $n$ 
with only real and distinct zeros. In the above notation the
following system of inequalities holds
         $$x^{(i)}_l<x^{(j)}_l<x^{(i)}_{l+j-i};\quad i<j\le n-l. \tag1$$
    We call this system
{\it the standard Rolle's restrictions.}

With any $n$-nice function $f$ we can associate the arrangement
$\Cal A_{f}$ of all $\binom {n+1} 2$ zeros $\{x^{(i)}_{l}\}$ of
$f^{(i)},\,i=0,\ldots,n-1;\,1\le l\le n-i $, say, taking first all 
$x^{(0)}_{l}$ then all $x^{(1)}_{l}$ etc. 
\medskip 
 
The main problem we address in this note is as follows. 
\medskip 
{\smc Question.} What additional restrictions besides $(1)$ exist 
on the arrangements $\Cal A_{f}=\{x^{(i)}_{l}\}$ for 
$n$-nice functions?  Or even more ambitiously, given an arrangement of 
$\binom {n+1} 2$ real numbers $\Cal A=\{x^{(i)}_{l}\vert\; 
i=0,\ldots,n-1;\,l=1,\ldots n-i\}$ satisfying the standard Rolle's 
restrictions  is it possible to say if there exist an $n$-nice function $f$ such that 
$\Cal A_{f}=\Cal A$? 
\medskip 

Our choice of the class of $n$-nice functions is motivated by the (easy to formalize) idea  that as soon as 
 one allows several real zeros of $f'$ in between two consecutive zeros of $f$ then no  interesting additional restrictions 
 are possible. Also this class seems to be the natural generalization of the well-studied class of real polynomials 
with all real zeros.
 
We will soon discover  that  the fact that a smooth function is 
$n$-nice implies additional inequalities on the components of 
$\Cal A_{f}$. Notice that the set of all $n$-nice functions on a 
given interval forms the subset of a larger class of functions 
which in several aspects behave similarly to real 
polynomials of degree $n$.   Let us now define this smooth analog of polynomials.  
We postpone the discussion of the usual polynomials 
until the next section. 
 
     \medskip 
     {\smc Main Definition.} A smooth real-valued function $f$ defined 
     on some interval $I$   is called {\it a pseudopolynomial of  degree $n$\/} 
     on $I$  if $f^{(n)}$ never vanishes. 
 
\medskip 
The usual Rolle's theorem  immediately implies that any 
 pseudopolynomial $f$ of degree $n$ has at most $n$ real zeros  (counted with 
multiplicitites). The set of all 
$n$-nice functions on a given interval coincides with the set of all 
 pseudopolynomials of degree $n$ with exactly $n$ real and distinct zeros. 
 \medskip
 
Let  $\Cal N_{n}(I)$ denote the set of all
   $n$-nice functions on an interval $I$. (Since the particular choice of $I$
   is unimportant in the formulation below we will often use $\Cal
   N_{n}$ instead of $\Cal N_{n}(I)$.)
    Our result below answers the posed question for $n=3$, i.e., the case of $3$-nice functions.  
   (Obviously, there are no additional restrictions for  
   $n=2$, i.e. on the arrangements of 
  the two zeros of a smooth function  $f$ and one zero of $f'$.) 
   In order to make our notation more readable 
   denote the three zeros of a $3$-nice function $f$  by 
$x_{1}<x_{2}<x_{3}$, the two zeros of $f'$ by $y_{1}< y_{2}$ and the only 
zero of $f''$ by $z_{1}$.

\medskip

{\smc Inequality Theorem. } For any $3$-nice function $f$ its  arrangement 
$\Cal A_{f}=(x_{1},x_{2},x_{3},y_{1},y_{2},z_{1})$ satisfies the 
following inequalities:

$$\cases \aligned 
              &x_{1}<y_{1}<x_{2}<y_{2}<x_{3}\\ 
              &y_{1}<z_{1}<y_{2}\\
              & y_{1}-x_{1}< \min\left(x_{2}-y_{1}, \sqrt{(z_1-y_1)^2+2(z_1-y_1)\vert z_1-x_2\vert}\right)\\ 
&x_{3}-y_{2}< \min\left( y_{2}-x_{2},\sqrt{(y_2-z_1)^2+2(y_2-z_1)\vert z_1-x_2\vert}\right)\\ 
\endaligned \endcases\tag2
 $$

And, conversely, for any $6$-tuple satisfying the above inequalities 
there exists a $3$-nice function $f$ with that 
arrangement $\Cal A_{f}$ of the zeros of $f,f',f''$ respectively. 
(The geometrical meaning of the additional inequalities will be quite 
  clear from the proof of the theorem below.) 
\medskip

The inequalities on the first two lines of $(2)$ are the standard Rolle's restrictions. Notice that two  remaining new inequalities interchange places under the  substitution $x\mapsto -x$. 

{\smc  Proof.}  Take some $f\in \Cal N_{3}$.  Without loss of generality we can 
assume  $f^{\prime\prime\prime}>0$  implying that $f'$ is 
convex (otherwise multiply $f$ by $-1$).  As above, the three 
zeros of $f$ 
are denoted 
by $x_{1}<x_{2}<x_{3}$, the two zeros of $f'$ by $y_{1}< y_{2}$ and the only 
zero of $f''$ by $z_{1}$. 
Let us first consider the case  $x_{2}<z_{1}$; see Fig.1.

    \medskip 
\vskip 15pt 
\centerline{\hbox{\epsfysize=5.5cm\epsfbox{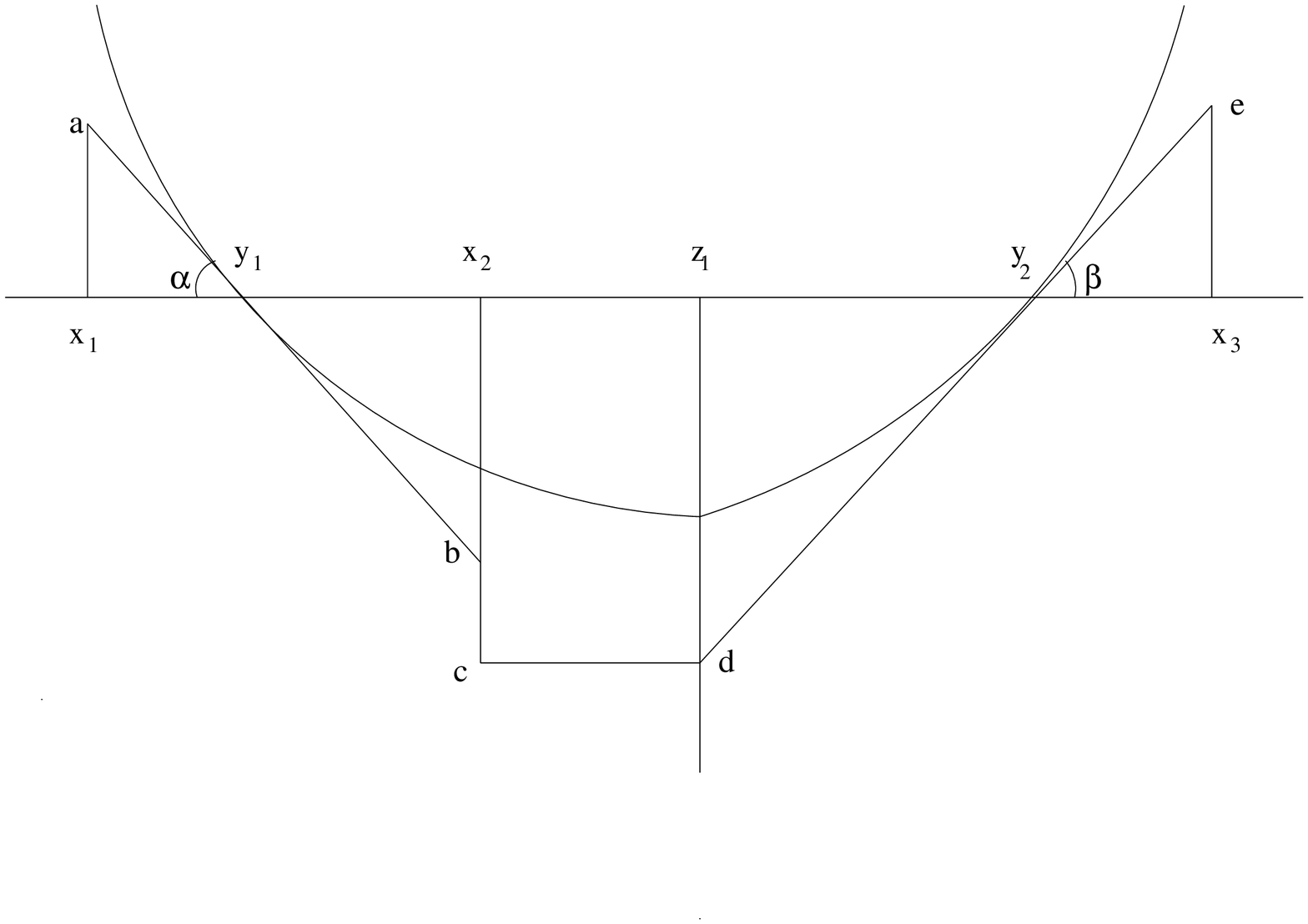}}} 
\midspace{0.1mm} \caption{Fig.~1. Derivative of $f$ and its accompanying 
    elementary configuration.} 
  
 \medskip
 
 We will
 immediately derive the additional inequalities: 
 $$ y_{1}-x_{1}< x_{2}-y_{1}\quad\text{ and }\quad x_{3}-y_{2} <\sqrt{(y_{2}-z_{1})^2+2(y_{2}-z_{1})(z_1-x_{2})}.\tag3$$ 
 The case $z_{1}<x_{2}$ is completely analogous and leads to: 
 $$x_{3}-y_{2}<y_{2}-x_{2}\quad\text{ and }\quad y_{1}-x_{1} <\sqrt{(z_{1}-y_{1})^2+2(z_1-y_1)(x_{2}-z_1)}.\tag4$$
    The union of $(3)$ and $(4)$ gives exactly the required new inequalities in  $(2)$. 
 
    Indeed, the convexity of $f'$ implies that  

\noindent 
    i)  on  $(x_{1},y_{1})$ the graph of $f'$  lies above the 
line segment $\{a,y_{1}\}$ which is  tangent to it at $y_{1}$; 
 
\noindent 
ii) on  $(y_{1},x_{2})$ the graph of $f'$ lies between 
the $x$-axis and the  line segment $\{y_{1},b\}$ which is tangent to it at $y_{1}$; 
 
\noindent 
iii) on  $(x_{2},y_{2})$ the graph of $f'$ lies between the 
$x$-axis and the broken line segment $\{c,d,y_{2}\}$ where $(c,d)$ is 
horizontal; 
 
\noindent 
iv)  on  $(y_{2},x_{3})$ the graph of $f'$  lies 
above the  line segment $\{y_{2},e\}$ which is tangent to it at $y_{2}$. 
 
Note that since $x_{1},x_{2},x_{3}$ are consecutive real zeros of $f$  then 
$$\int_{x_1}^{x_2}f'dx=\int_{x_2}^{x_3}f'dx=0$$ 
implying  
$$\int_{x_{1}}^{y_{1}}f'dx=-\int_{y_{1}}^{x_{2}}f'dx\quad \text{ and }\quad  \int_{x_{2}}^{y_{2}}f'dx=-\int_{y_{2}}^{x_{3}}f'dx.$$ 
    Therefore, 
    $$Ar(\triangle_{ax_{1}y_{1}}) < Ar(\triangle_{y_{1}bx_{2}})\quad 
   \text{ and }\quad  Ar(\triangle_{ey_{2}x_{3}}) < Ar(\square_{x_{2}cdy_{2}}),$$
    where 
    $Ar$ stands for the area of the corresponding figures. (The figure 
    $\square_{x_{2}cdy_{2}}$ is a trapezoid.) Using our 
    knowledge of high-school mathematics we get 
    $$Ar(\triangle_{ax_{1}y_{1}})=\frac {1}{2}(y_{1}-x_{1})^2\tan \alpha,\; 
    Ar(\triangle_{ey_{2}x_{3}})=\frac {1}{2}(x_{3}-y_{2})^2\tan \beta,\; 
    Ar(\triangle_{y_{1}bx_{2}})=\frac{1}{2}(x_{2}-y_{2})^2\tan \alpha, $$ 
    and,  
    $$Ar(\square_{x_{2}cdy_{2}})=Ar(\triangle_{z_{1}dy_{2}}) 
    +Ar(\square_{x_{2}cdz_{1}})=\frac{1}{2}(y_{2}-z_{1})^2\tan 
    \beta+(y_{2}-z_{1})(z_{1}-x_{2})\tan \beta=$$ $$=\frac {1}{2}\left((y_{2}-z_{1})^2+2(y_{2}-z_{1})(z_{1}-x_{2})\right)\tan \beta.$$
     These relations  immediately imply the required inequalities. 
 
   To finish the proof, pick any $6$-tuple of real numbers 
    $x_{1},x_{2},x_{3},y_{1},y_{2},z_{1}$ satisfying $(2)$. We 
    again assume that $x_{2}< z_{1}$. (The case $z_{1}\le x_{2}$ is completely 
    analogous.) Draw a piecewise linear function $g$ as shown on Fig. 2.

    \medskip 
\vskip 15pt 
\centerline{\hbox{\epsfysize=3.0cm\epsfbox{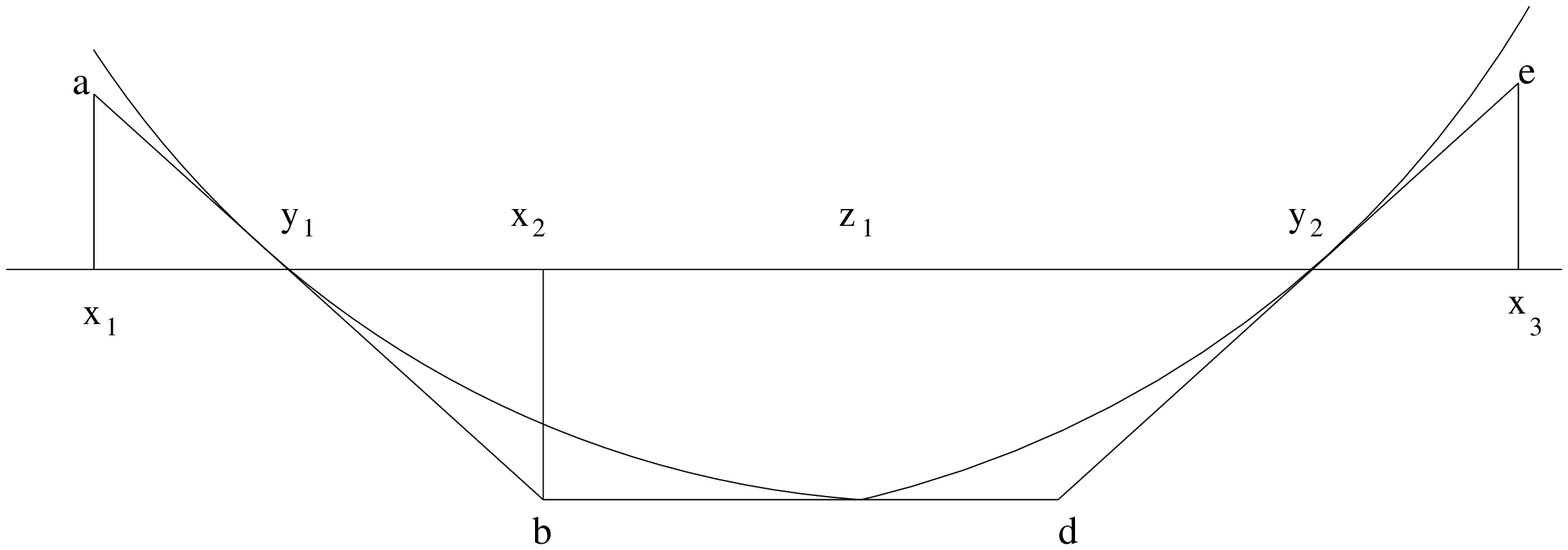}}} 
\midspace{0.1mm} \caption{Fig.~2. Constructing an appropriate function.} 
\medskip 
    (The only difference with Fig.1 is that we force $c=b$.) 
    The inequalities $(2)$ imply that $Ar(\triangle_{ax_{1}y_{1}}) < 
Ar(\triangle_{y_{1}bx_{2}})$ 
    and $Ar(\triangle_{ey_{2}x_{3}}) < Ar(\square_{x_{2}bdy_{2}})$. 
    It is easy to approximate $g$ by a convex function $h$ so that 
     
     i) $y_{1}$ and $y_{2})$ are the zeros of $h$ and $z_{1}$ is the zero
     of $h'$;
 
     ii) $\int_{x_{1}}^{y_{1}}hdx<-\int_{y_{1}}^{x_{2}}hdx$ and
 $\int_{x_{2}}^{y_{2}}hdx<-\int_{y_{2}}^{x_{3}}hdx$.
 
 Keeping the function $h$ as it is on the interval $(y_{1},y_{2})$
 we can increase it  on the intervals $(x_{1},y_{1})$ and
 $(y_{2},x_{3})$ and construct  a new convex function $\tilde h$
  with the properties 
 
     i) $y_{1}$ and $y_{2}$ are the zeros of $\tilde h$ and $z_{1}$
 is the zero
     of $\tilde h'$;
 
      ii) $\int_{x_{1}}^{y_{1}}\tilde hdx=-\int_{y_{1}}^{x_{2}}\tilde
 hdx$ and
 $\int_{x_{2}}^{y_{2}}\tilde hdx=-\int_{y_{2}}^{x_{3}}\tilde hdx$.
 
 The function $f\in \Cal N_{3}$ we were looking for is now obtained
 as $f=-\int_{y_{1}}^{x}\tilde h(t)dt$. \qed
 \medskip 
 
 As an illustration of the Inequality Theorem consider the polynomial 
 $$p(x)=x(x-1)(x-4).$$
 The polynomial $p(x)$ and its first two derivatives are shown on Fig.3 below. 
 As any polynomial of degree $3$ with real and distinct zeros $p(x)$ is $3$-nice. Elementary calculation gives that  in the notation of the theorem we have 
 $$x_1=0, x_2=1, x_3=4,\; y_1=\frac{5-\sqrt{13}}{3}, y_2=\frac{5+\sqrt{13}}{3},\; z_1=\frac  5 3.$$ 
 Therefore, $x_2<z_1$ and we are interested in checking the validity of $(3)$. 
 Indeed, $y_1-x_1=\frac{5-\sqrt{13}}{3}\simeq 0.464816$ which is smaller than $x_2-y_1=\frac{\sqrt{13}-2}{3}\simeq 0.535184$.  Next, $x_3-y_2=\frac{7-\sqrt{13}}{3}\simeq 1.13148$ which is smaller than $\sqrt{(z_{1}-y_{1})^2+2(z_1-y_1)(x_{2}-z_1)}= \frac{\sqrt{13+4\sqrt{13}}}{3}\simeq 1.74554.$

  \medskip 
\vskip 15pt 
\centerline{\hbox{\epsfysize=3.5cm\epsfbox{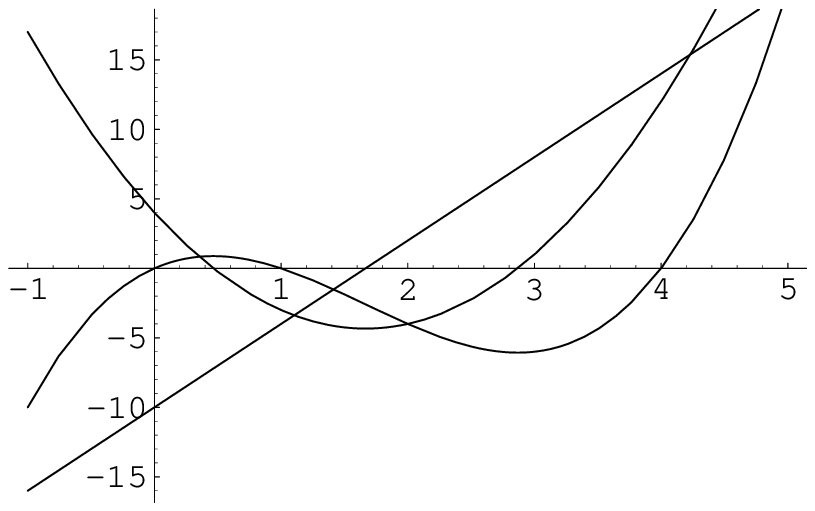}}} 
\midspace{0.1mm} \caption{Fig.~3. Illustration of the theorem.} 
\medskip 
 
      \heading  Usual polynomials \endheading

    Let us now try to find meaningful restrictions on arrangements ${\Cal A}_{p}$ 
     for the usual polynomials with all real and distinct zeros. 
     (Note that inequalities $(2)$ remain valid but they clearly do not give all the restrictions 
     since for a usual polynomial its zeros define 
     all the zeros of all its higher derivatives in a unique way.) 
     Let us assume that all $n$ zeros $x_{1}^{(0)}<x_{2}^{(0)}<\ldots<x_{n}^{(0)}$ 
     of a  polynomial $p$ of degree $n$ are real and distinct. 
     Assume additionally that all $x_{l}^{(i)}$ are pairwise different. 
     (One can see 
     that this extra condition holds for almost all polynomials with 
     real zeros.)  We  call such polynomials {\it strictly $n$-nice}.  For a 
     strictly $n$-nice polynomial $p$ 
     its whole arrangement $\Cal A_{p}$ is naturally ordered on the real line. 
     Substituting each zero of $p$ by the symbol $0$, each zero 
     of $p'$ by $1$, ... , each zero of $p^{(n-1)}$ by $(n-1)$ respectively we get a {\it symbolic sequence of $p$} of 
length $\binom {n+1} {2}$ 
     with $n$ occurrences of $0$, $(n-1)$ occurrences of $1$,..., one occurence of $(n-1)$  
     and satisfying the condition that between any two consecutive occurrences 
     of the symbol $i$ it has exactly one occurrence of the symbol $i+1$. 
 
     For example, there are only two possible symbolic sequences  for 
     $n=3$, namely, $012010$ and $010210$. For $n=4$ there are 12 such 
     sequences $0123012010$, $0120312010$, $0120132010$, $0102312010$, 
     $0102132010$, 
     $0123010210$, $0120310210$, $0120130210$, $0120103210$, $0123010210$,       
     $0102130210$, $0102103210$. A patient reader will find the for $n=5$ there are 
     $286$ such sequences. 
 
  If we denote by $\flat_{n}$ the number of all possible symbolic sequences 
  of length $n$ then 
actually, this number is possible to calculate. 
It turns out to be equal to 
$$\flat_n=\binom{n+1} 
{2}!\frac{1!2!...(n-1)!}{1!3!...(2n-1)!}.$$  
But since this calculation is a content of a different story we refer the interested reader to \cite  {6}. 
 
We can now formulate  a natural discrete analog of the main problem from the previous section which makes sense for the usual polynomials.  
\medskip 
{\smc Question.} What symbolic sequences can occur for strictly 
$n$-nice polynomials of degree $n$? (We will call such sequences 
{\it realizable}.) 
\medskip 
 
As we will see shortly already the first nontrivial case of $n=4$ shows that 
the number $\sharp_{n}$ of all realizable symbolic sequences is 
strictly smaller than the corresponding $\flat_{n}$, namely, 
$10=\sharp_4\neq\flat_4=12$.

\medskip 
The fact that $\sharp_n\neq\flat_n$ 
was apparently observed by a number of authors but the only relevant 
reference we found is \cite {1} published in 1993. An 
explanation of this phenomenon for $n=4$ is as follows. 
(See further generalizations in \cite {4}.) 
\medskip 
{\smc Theorem,} see  \cite {1}. A  polynomial $p$ of degree $4$ 
with real zeros $x_{1}<x_{2}<x_{3}<x_{4}$ satisfying 
the inequalities $x_{2}<z_{1}$ and $x_{3}<z_{2},$  satisfies additionally the 
inequality $y_{2}<t_{1}$. Here $y_{1}<y_{2}<y_{3}$ are the zeros of $f'$; 
  the zeros of $f''$ are $z_{1}<z_{2}$  and $t_{1}$ is the zero of 
$f^{\prime\prime\prime}$. In other words, the symbolic sequences 
$0102310210$ and $0120132010$ are non-realizable. 
\medskip 
 
This is easy to check once you know what to prove! Indeed,  any monic 
polynomial of degree $4$ with all real zeros can be put in the 
form $x^4-x^2+ux+v$ by a linear change of $x$ and scaling. Namely, 
by shifting $x\mapsto x+\al$ we can always get rid of the 
$x^3$-term. Since the second derivative of the  obtained polynomial 
 has two real zeros, the coefficient at $x^2$ should be 
negative. Appropriate scaling now puts $p$ in the above form. Note that 
$p''=12x^2-2$, its zeros being $\pm \sqrt{\frac{1}{6}}$. The 
assumptions $x_{2}<z_{1},\;x_{3}<z_{2}$ together with $p$ having 
real zeros imply $p(- \sqrt{\frac{1}{6}})>0,\; p( 
\sqrt{\frac{1}{6}})<0$ (draw the graph of $p$). Noting that 
$p^{\prime\prime\prime}=24x$ what we need to prove is that 
$p'(0)<0$ (draw the graph of $p'$). The last inequality is 
equivalent to $u<0$. Expanding $p(- \sqrt{\frac{1}{6}})>0$ and 
$p(\sqrt{\frac{1}{6}})<0$ we get $\frac {1}{36}-\frac {1}{6}-\frac 
{u}{\sqrt{6}}+v>0$ and $\frac {1}{36}-\frac {1}{6}+\frac 
{u}{\sqrt{6}}+v<0$. Subtracting the former from the later implies 
$\frac {2u}{\sqrt 6}<0$. \qed 
 
\medskip 
 
The next case $n=5$ was considered in \cite {2}. V.~Kostov was able to show that among $286$ possible symbolic sequences only $116$ are realizable by strictly $n$-nice polynomials. Very recently the same author considered the similar question  which symbolic sequences are realizable for the case of $n$-nice functions; see \cite {3}. It turned out that for $n=4$ all $12$ symbolic sequences are realizable but already for $n=5$ there are non-realizable sequences. The situation  does not seem to change much if we extend the class of polynomials but $n$-nice functions. 

To finish the section let us present a tempting problem  posed by the famous mathematician Vladimir Arnold after the talk given by V.~Kostov on his seminar. 

\medskip 
{\smc Problem.} Is it true that $\lim_{n\to\infty} \frac{\sharp_n}{\flat_n}=0$? If yes, how fast does the quotient 
$\frac{\sharp_n}{\flat_n}$ decrease?  
\medskip 
 
 \heading Periodic functions \endheading
      
     At the end  let us  briefly discuss what happens with periodic functions, i.e., functions defined on a  circle. 
      In the previous sections we defined the class of $n$-nice functions --
      a generalization of polynomial of degree $n$
      with $n$ distinct real roots -- and
      found some inequalities involving the roots of higher derivatives  
      valid    for any $3$-nice smooth function.  It seems quite natural to try to develop a similar concept for periodic  
       functions.  The  periodic analog of  polynomials of degree $n$ are 
       trigonometric polynomials of degree $n$, i.e., expressions of the 
      form $a_{0}+\sum_{k=1}^n a_{k}\cos kx + b_{k}\sin kx$. Any 
      trigonometric polynomial of  degree $n$ has at most $2n$ real zeros 
      on a period. 
      
      Observe that if we take such a trigonometric polynomial 
      with exactly $2n$ real and distinct zeros then  its derivative of {\sl any} order will also be a trigonometric polynomial of the same degree $n$. Moreover, by the usual   Rolle's theorem it will also have  exactly  $2n$ real and 
      distinct zeros on a period. So it seems tempting to define a 
      periodic analog of an $n$-nice function as a periodic function 
      such that it and its derivatives of any order have exactly $2n$  
      real zeros. But (for not completely clear reasons) the situation with periodic functions 
      turns out to be much more rigid than with the functions on an interval.  
      
      To explain the situation we have to invoke the following famous classical result of G.~Polya and N.~Wiener; see 
      \cite {5}. 
     \medskip 
     {\smc Theorem.} Any periodic function $f$ such that the number of 
     real  zeros of the $i$-th derivative $f^{(i)}$ on a period remains 
bounded as $i\to \infty$ is   a trigonometric polynomial. 
   \medskip 
 
 In particular, any (conjectural) $n$-nice periodic function must necessarily be an actual  trigonometric polynomial.  
 On the other hand, one can define a  periodic analog of  symbolic sequences from the previous section and ask which of those are realizable by trigonometric polynomials with all real zeros.  Namely, for a positive integer $k$ consider a sequence of integers of length $nk$  written on a circle and containing $n$ zeros, $n$ ones, $n$ twos,..., $n$ copies of $(k-1)$. We call such a sequence {\it possible periodic } if for any $i=0,1,...,k-2$ in between 
 any two consecutive (on the circle) copies of $i$ the sequence contains exactly one copy of $i+1$. Now  we can
ask: 

\medskip 
{\smc Question.} What possible periodic sequences can occur as the sequences of zeros of $f$ and its higher derivatives  of order up to $k-1$ where $f$ is a trigonometric polynomial of degree $n$ with all real and distinct zeros? (Here as before  the integer $i$ substitutes  a real root of the $i$-th derivative of $f$.) 
\medskip 

Unfortunately, at the moment there is no nontrivial information available about the latter problem.   
 \medskip 
 \noindent 
 {\smc Acknowledments.}Ê The authors are sincerely grateful to the anonymous referees whose important suggestions allowed us to substantially improve the quality of exposition.  
 

\bigskip 
 
\Refs 
\widestnumber \key{ShSh} 
 
\ref \key 1 \by B.~Anderson 
\paper {\rm Polynomial root dragging} 
\jour {\sl Amer. Math. Monthly} 
\vol 100 
\yr 1993 
\pages 864--866 
\endref 
 
\ref \key 2 \by V.Kostov 
\paper {\rm Discriminant sets of families of hyperbolic polynomials of degree $4$ and $5$ } 
\jour  {\sl Serdica Math. J} 
\vol 28 
\yr 2002 
\issue 2 
\pages 117--152 
\endref 

\ref \key 3 \by V.Kostov 
\paper {\rm On  polynomial-like functions} 
\jour  {\sl Bulletin des sciences mathematiques } 
\finalinfo to appear 
\endref 

 
\ref \key 4 \by V.~Kostov and B.~Shapiro 
\paper {\rm On arrangements of roots for a real hyperbolic polynomial and 
its derivatives} 
\jour {\sl Bulletin des sciences mathematiques }
\yr 2002 
\vol 126 
\issue 1 
\pages 45--60 
\endref 
 
\ref \key 5 \by G.~Polya and N.~Wiener 
\paper {\rm On the oscillation of the derivatives of a periodic function} 
\jour {\sl Trans. Amer. Math. Soc} 
\vol 52 
\yr 1942 
\pages 249--256 
\endref 
 
 
 
\ref \key 6 \by M.~Thrall 
\paper {\rm A combinatorial problem}Ê 
\jour {\sl Michigan Math. J} 
\vol 1 
\yr 1952 
\pages 81--88 
\endref 
 
  

\endRefs
\enddocument